\documentclass[a4paper,12pt,twoside,onecolumn]{article}
\usepackage{fix-cm}
\usepackage[OT1]{fontenc}

\usepackage{amsmath,amssymb,amsfonts,amsthm}
\usepackage{graphicx,float,accents}
\usepackage[a4paper,textwidth=16cm,textheight=24cm,centering]{geometry}
\usepackage[pagewise]{lineno}
\usepackage{fancyhdr,caption,titlesec,setspace,indentfirst}
\usepackage[style]{abstract}

\renewcommand{\abstitlestyle}[1]{}
\usepackage{etoolbox}
\makeatletter
\patchcmd{\@maketitle}{\par\vskip 1.5em}{\par\vskip 6pt}{}
  {\PackageError{paper-layout}{Unable to set title-block spacing}{Check the article title definition.}}
\makeatother
\usepackage[numbers,sort]{natbib}
\usepackage{enumitem}
\setlist[enumerate]{label=\textnormal{(\arabic*)},leftmargin=3em,
  labelsep=0.5em,itemsep=0pt,partopsep=0pt,parsep=0pt,topsep=3pt}
\setitemize[1]{itemsep=-4pt,partopsep=0pt,parsep=\parskip,topsep=5pt}
\setdescription{itemsep=-4pt,partopsep=0pt,parsep=\parskip,topsep=5pt}
\usepackage{microtype}
\usepackage[hyperfootnotes=false,colorlinks=true,linkcolor=blue,citecolor=blue,urlcolor=blue]{hyperref}
\allowdisplaybreaks
\titleformat{\section}[hang]
  {\normalfont\rmfamily\fontsize{14}{17}\selectfont\bfseries\raggedright}
  {\thesection}{1em}{}
\titlespacing{\section}{0pt}{\baselineskip}{0.5\baselineskip}
\titleformat{\subsection}[hang]
  {\normalfont\rmfamily\fontsize{12}{14.4}\selectfont\bfseries\raggedright}
  {\thesubsection}{1em}{}
\titlespacing{\subsection}{0pt}{\baselineskip}{0.5\baselineskip}
\titleformat{\subsubsection}[hang]
  {\normalfont\rmfamily\fontsize{12}{14.4}\selectfont\bfseries\raggedright}
  {\thesubsubsection}{1em}{}
\titlespacing{\subsubsection}{0pt}{\baselineskip}{0.5\baselineskip}
\titleformat{\paragraph}[hang]
  {\normalfont\rmfamily\fontsize{12}{14.4}\selectfont\bfseries\raggedright}
  {\theparagraph}{1em}{}
\titlespacing{\paragraph}{0pt}{\baselineskip}{0.5\baselineskip}
\titleformat{\subparagraph}[hang]
  {\normalfont\rmfamily\fontsize{12}{14.4}\selectfont\bfseries\raggedright}
  {\thesubparagraph}{1em}{}
\titlespacing{\subparagraph}{0pt}{\baselineskip}{0.5\baselineskip}
\newcommand{\footnoteONLYtext}[1]{\begingroup
\renewcommand{\thefootnote}{}\footnotetext{#1}\endgroup}
\newcommand{\keywords}[1]{{\it Keywords}: #1}
\newcommand{\msc}[1]{{\it 2020 Mathematics Subject Classification}: #1}
\numberwithin{equation}{section}
\theoremstyle{plain}
\newtheorem{theoalph}{Theorem}

\newtheorem{prop}{Proposition}[section]
\newtheorem{lemm}[prop]{Lemma}

\theoremstyle{definition}
\newtheorem{defi}[prop]{Definition}
\theoremstyle{remark}
\newtheorem{rema}[prop]{Remark}

\newcommand{\R}{\mathbb R}
\newcommand{\N}{\mathbb N}
\newcommand{\Nzero}{\mathbb N_0}
\newcommand{\aseq}{\mathbf a}
\newcommand{\T}{\mathcal T}
\newcommand{\dist}{\operatorname{dist}}
\newcommand{\Ea}{E_{\aseq}}
\newcommand{\Ta}{T_{\aseq}}
\newcommand{\card}{\#}
\hypersetup{pdftitle={Schmidt games for nonautonomous Liverani--Saussol--Vaienti maps}}
\title{\Large\bfseries Schmidt games for nonautonomous\\Liverani--Saussol--Vaienti maps}
\author{Haoyang Ji}
\date{}
\begin{document}
\maketitle
\begin{abstract}
\noindent\textbf{Abstract.} Liverani--Saussol--Vaienti (LSV) maps form
a family of piecewise differentiable dynamical systems on $[0,1]$
depending on a parameter $a>0$. These maps are expanding away from
a neutral fixed point. In this paper, we study Schmidt games for
nonautonomous compositions of LSV maps with arbitrary parameter
sequences in a compact subinterval of $(0,\infty)$.
For every countable set $P\subset[0,1]$, we prove that the set of
points whose orbit closure is disjoint from $P$ is absolute winning
on $[0,1]$. Our proof combines a potential game argument with an
interval tree constructed from successive returns.
\end{abstract}
\footnoteONLYtext{Date: 30 October 2026}
\footnoteONLYtext{\msc{37E05, 91A44}}
\footnoteONLYtext{\keywords{intermittency, LSV maps, nonautonomous dynamics, absolute winning, Schmidt games}}
\footnoteONLYtext{H. Ji was supported by NSFC Grant No.12301103.}
\section{Introduction}

Let us consider the one-parameter family of LSV maps. For $a >0$, let 
\begin{equation}\label{eq:lsv}
T_a(x)=
\begin{cases}
x(1+2^a x^a),&0\le x<1/2,\\
2x-1,&1/2\le x\le1.
\end{cases}
\end{equation}
These maps were introduced by Liverani, Saussol, and Vaienti \cite{LSV}
as a simplified version of the Pomeau--Manneville family \cite{PM80}.
Since then they have become a standard model of intermittency.
They have two monotonicity intervals independent of $a$, a neutral
fixed point at $0$, and an affine right branch of slope $2$.
Orbits may spend long periods near the neutral fixed point at $0$,
where expansion is weak, while the right branch expands by a factor of $2$. 
These maps have been studied extensively. For decay of correlations,
see \cite{LSV,Young,Hu}; for linear response, see \cite{BT,Korepanov};
and for large deviations, see \cite{Mel}, among others.
For stochastic stability of related intermittent circle maps,
see \cite{AT,SvS}; for statistical properties of random LSV
compositions, see \cite{BBD,BBR,BQT}.

In this paper we mainly study Schmidt games for nonautonomous compositions of
LSV maps. 
Fix
\[
0<a_-\le a_+<\infty, \ 
K=[a_-,a_+], \  \aseq=(a_n)_{n\ge0}\in K^{\Nzero}.
\]
Set
\[
\Ta^0=\mathrm{id}, \ 
\Ta^n=T_{a_{n-1}}\circ\cdots\circ T_{a_0} \ (n\ge1).
\]
For $p\in[0,1]$ and $P\subset[0,1]$, define the exceptional sets
\begin{align}
\Ea(p)&=\{x\in[0,1]:\inf_{n\ge0}|\Ta^n(x)-p|>0\},\label{eq:single}\\
\Ea(P)&=\bigcap_{p\in P}\Ea(p).\label{eq:countable}
\end{align}
Equivalently,
\[
\Ea(P)=
\{x\in[0,1]:
P\cap\overline{\{\Ta^n(x):n\ge0\}}=\varnothing\}.
\]
The main theorem of this paper is the following.

\begin{theoalph}\label{thm:main}
Fix $0<a_-\le a_+<\infty$ and $K=[a_-,a_+]$.
For each sequence $\aseq\in K^{\Nzero}$ and each countable
set $P\subset[0,1]$, the set $\Ea(P)$ is absolute winning on
$[0,1]$.
\end{theoalph}
In particular, $\Ea(P)$ is strong winning and Schmidt winning on $[0,1]$, and has full
Hausdorff dimension in every nondegenerate subinterval of $[0,1]$;
see \cite{McM,Schmidt}.

Schmidt's game \cite{Schmidt} and its variants \cite{McM,FPS10} provide a method for
studying the size of sets of points with nondense orbits. For expanding dynamical systems, see, for example, \cite{Tseng,Farm,FPS10,Wu,HLY}.
For intermittent interval maps, Duvall \cite{Duv20} proved strong
winning of the exceptional set for the neutral fixed point.
By applying Hu--Li--Yu's theorem \cite{HLY} and an approximation
argument, Duvall established strong winning for arbitrary targets
\cite{DuvJun} and absolute winning for countable target sets
\cite{DuvJul}. For winning properties of exceptional sets in certain
nonautonomous linear systems, see \cite{BBFKW,BFK}.
Huang--Li--Wang--Yuan \cite{HLWY} studied a class of nonautonomous
product maps on $[0,1]^d$, whose one-dimensional case consists of
nonautonomous piecewise expanding interval maps.

Our system in considerable is not uniformly expanding, so the previous results, in particular \cite[Theorem 1.4]{HLWY}, do not
apply directly.
Our proof follows the potential game approach and inverse branch
counting method of Wu \cite{Wu}.
We construct an interval tree by pulling back $Y=[1/2,1]$ along successive
return branches. The tree structure provides the counting estimate needed
to construct Alice's winning strategy.
Unlike Duvall's proofs \cite{DuvJun,DuvJul}, our argument does not
transfer a winning strategy from an induced map to the original system.
For the standard LSV maps, our argument
provides a shorter and different proof of Duvall's absolute-winning
result \cite{DuvJul}.

This article is organized as follows.
Section~\ref{sec:games} establishes an absolute-winning criterion
for marked interval trees.
Section~\ref{sec:lsv} constructs the return interval tree and
proves
uniform distortion bounds for its branch maps and uniform
contraction of child intervals.
Section~\ref{sec:avoidance} proves Theorem~\ref{thm:main},
first treating the neutral fixed point and then countable
target sets.

\section{Absolute winning on interval trees}\label{sec:games}
\subsection{The games}

We first recall McMullen's absolute game \cite{McM}, using the
convention of \cite{BHNS}.
We write $\overline{B}(x,r)=[x-r,x+r]$, so that $|\overline{B}(x,r)|=2r$, where $|\cdot|$ denotes
length (Lebesgue measure).

\begin{defi}[Absolute game]\label{def:game}
Let $I$ be either $\R$ or a nondegenerate compact interval,
and let $S\subset I$.
Fix $0<\beta<1/3$. In the $\beta$-absolute game, Bob chooses a
closed interval $B_1=\overline{B}(x_1,r_1)\subset I$ with $r_1>0$.
After Bob chooses $B_n=\overline{B}(x_n,r_n)$, Alice deletes a closed interval
$A_n=\overline{B}(y_n,\rho_n)$ satisfying
\[
0<\rho_n\le\beta r_n.
\]
Bob then chooses a closed interval $B_{n+1}=\overline{B}(x_{n+1},r_{n+1})$
satisfying
\[
B_{n+1}\subset B_n\setminus A_n, \ 
r_{n+1}\ge\beta r_n.
\]
Alice wins by default if $r_n\not\to0$.
Otherwise, $\bigcap_{n\ge1}B_n=\{x_\infty\}$, and Alice wins
if $x_\infty\in S$.
If Alice can win regardless of Bob's choices, $S$ is called
{\it $\beta$-absolute winning} on $I$. It is called {\it absolute winning} on $I$
if it is $\beta$-absolute winning on $I$ for every $0<\beta<1/3$.
For $I=\R$, this defines absolute winning on $\R$.
We omit the ambient interval when it is clear from the context.
\end{defi}

We recall the standard properties of absolute winning sets;
see \cite[Section~1]{McM}.

\begin{prop}\label{prop:aw-properties}
Let $I$ be either $\R$ or a nondegenerate compact interval.
\begin{enumerate}
\item If $S\subset I$ is absolute winning and
$S\subset S'\subset I$, then $S'$ is absolute winning.
\item If $S_k\subset I$ is absolute winning for every
$k\in\N$, then $\bigcap_{k\in\N}S_k$ is absolute winning.
\item If $S\subset I$ is absolute winning and
$\mathcal C\subset I$ is countable, then $S\setminus\mathcal C$
is absolute winning.
\item Every absolute winning subset of $I$ is strong winning and
Schmidt winning.
\end{enumerate}
\end{prop}

We use the point version of the potential game introduced in
\cite[Appendix~C]{FSU}; see also \cite[Section~2.3]{BHNS}.

\begin{defi}[Potential game]\label{def:potential-game}
Let $S\subset\R$. Fix $s>0$ and $\beta\in(0,1)$.
In the $(s,\beta)$-potential game on $\R$, Bob first chooses
a closed interval $B_1=\overline{B}(x_1,r_1)\subset\R$
with $r_1>0$.
After $B_n=\overline{B}(x_n,r_n)$ is chosen, Alice specifies an at most
countable collection of closed intervals
$A_{n,j}=\overline{B}(y_{n,j},\rho_{n,j})$, of positive radii, such that
\begin{equation}\label{eq:budget}
\sum_j\rho_{n,j}^{s}\le(\beta r_n)^s.
\end{equation}
The collection may be empty. Bob then chooses a closed interval
$B_{n+1}=\overline{B}(x_{n+1},r_{n+1})$ satisfying
\[
B_{n+1}\subset B_n, \ 
r_{n+1}\ge\beta r_n.
\]
He need not avoid any of Alice's intervals.
Alice wins by default if $r_n\not\to0$.
Otherwise, $\bigcap_{n\ge1}B_n=\{x_\infty\}$, and Alice wins
for $S$ precisely when
\[
x_\infty\in S\cup\bigcup_{n\ge1}\bigcup_j A_{n,j}.
\]
The set $S$ is called {\it $(s,\beta)$-potential winning} if Alice has
a winning strategy. It is called {\it $0$-potential winning} if it is
$(s,\beta)$-potential winning for every $s>0$ and every
$\beta\in(0,1)$.
\end{defi}

The following proposition is due to Fishman, Simmons, and
Urba\'nski \cite[Appendix~C, Observation~C.7 and Theorem~C.8]{FSU};
see also \cite[Proposition~1.3]{BHNS}.

\begin{prop}\label{prop:conversion}
For $S\subset\R$,
\[
S\text{ is absolute winning}
\ \Longleftrightarrow\ 
S\text{ is }0\text{-potential winning}.
\]
\end{prop}

In particular, to prove that a set is absolute winning,
it suffices to prove that it is $(s,\beta)$-potential winning
for every $s>0$ and $\beta\in(0,1)$.
We use this criterion in the next subsection.

\subsection{Marked interval trees}\label{sec:trees}

We first recall the terminology for rooted trees; see
\cite{ABV}. Tree constructions have been used in \cite{ABV}
and \cite{BH}. Our definition differs in allowing
a vertex to have countably infinitely many children.

\begin{defi}[Rooted tree]\label{def:rooted-tree}
A rooted tree is a connected graph $\T$ without cycles, with a
specified vertex $\tau_0$ called its {\it root}. We identify $\T$
with its vertex set. The {\it depth} of $\tau\in\T$ is the number of
edges in the unique path from $\tau_0$ to $\tau$.
Let $\T_n$ denote the vertices of depth $n$. Thus
\[
\T=\bigcup_{n\ge0}\T_n, \  \T_0=\{\tau_0\}.
\]
We write $\tau\preceq\sigma$ if $\tau$ lies on the path from
$\tau_0$ to $\sigma$. In this case, $\tau$ is an {\it ancestor} of
$\sigma$, and $\sigma$ is a {\it descendant} of $\tau$; both relations
include equality. For $\tau\in\T_n$, its {\it children} form the set
\[
\T_{\mathrm{ch}}(\tau)
=\{\sigma\in\T_{n+1}:\tau\preceq\sigma\}.
\]
\end{defi}

\begin{defi}[Interval tree]\label{def:tree}
An {\it interval tree} on a nondegenerate compact interval $I\subset\R$
consists of a countable rooted tree $\T$ with root $\tau_0$ and
a family $(J_\tau)_{\tau\in\T}$ of nondegenerate closed subintervals
of $I$, with $J_{\tau_0}=I$, satisfying the following conditions:
\[
\begin{aligned}
\tau\preceq\sigma
&\ \Longrightarrow\ J_\sigma\subset J_\tau,\\
\tau\npreceq\sigma\ \text{and}\ \sigma\npreceq\tau
&\ \Longrightarrow\ J_\tau^\circ\cap J_\sigma^\circ=\varnothing,
\end{aligned}
\]
and, for some $q\in(0,1)$,
\begin{equation}\label{eq:contraction}
|J_\sigma|\le q|J_\tau|
 \ \mbox{ for all } \tau\in\T,\ \sigma\in\T_{\mathrm{ch}}(\tau).
\end{equation}
We identify each vertex $\tau$ with its interval $J_\tau$
and write $J\in\T$ for an interval of the tree.
Each vertex may have countably infinitely many children.

An interval tree is called {\it marked} if each interval $J\in\T$
is assigned a point $z_J\in J$.
\end{defi}

The following proposition provides the main criterion for proving
the absolute winning property.

\begin{prop}\label{prop:tree}
Let $\T$ be a marked interval tree with root interval $I$, let $z_J\in J$ be the marked points assigned to each $J\in\T$. Then the set
\[
\widehat S=
\left\{x\in\R:
\inf_{J\in\T}\frac{|x-z_J|}{|J|}>0\right\}
\]
is absolute winning on $\R$. In particular,
$\widehat S\cap I^\circ\ne\varnothing$.
\end{prop}

\begin{proof}
By Proposition~\ref{prop:conversion}, it suffices to show that
$\widehat S$ is $0$-potential winning.
Fix arbitrary $s>0$ and $\beta\in(0,1)$. Choose an integer $M\ge1$ with
$q^M<\beta$. Set $C_\beta=6M$, and choose
\begin{equation}\label{eq:delta}
0<\delta\le
\min\{\beta,\ \beta^2(2C_\beta)^{-1/s}\}.
\end{equation}
Since $\T$ is countable, choose a sequence $(J_n)_{n\ge1}$
containing every vertex of $\T$.

Let $B_1=\overline{B}(x_1,r_1)$ be Bob's first choice.
At each turn $n\ge1$, when Bob chooses $B_n=\overline{B}(x_n,r_n)$,
Alice first chooses the ball
$\overline{B}(z_{J_n},2^{-1/s}\beta r_n)$. She then chooses all balls
$\overline{B}(z_J,\delta|J|)$ satisfying
\begin{equation}\label{eq:band}
r_n\le|J|<r_n/\beta, \ 
\overline{B}(z_J,\delta|J|)\cap B_n\ne\varnothing.
\end{equation}

We claim that the number of balls satisfying \eqref{eq:band}
is bounded by $C_\beta$. Partition $\T$ into $M$ groups:
\[
\T=\bigcup_{i=1}^{M}\widehat{\T}_i, \ 
\widehat{\T}_i=\bigcup_{k\ge0}\T_{i-1+kM}
\  (1\le i\le M).
\]
Fix $1 \leq i\leq M$. Suppose that distinct vertices $J,J'\in\widehat{\T}_i$
satisfy \eqref{eq:band} and that $J'$ is a descendant of $J$. Then their depths differ by a positive multiple of $M$, so
\[
|J'|\le q^M|J|<\beta(r_n/\beta)=r_n,
\]
contradicting \eqref{eq:band}. Thus the intervals in
$\widehat{\T}_i$ satisfying \eqref{eq:band} have pairwise
disjoint interiors. Since $\delta\le\beta$, by \eqref{eq:band} we have
\[
|z_J-x_n|\le r_n+\delta|J|<2r_n.
\]
Then each such $J$ intersects
$W=[x_n-2r_n,x_n+2r_n]$ and has length at least $r_n$. There are at most 6 such intervals in each group $\widehat{\T}_i$, and hence the total number is at most $6M$.

At each turn $n\ge1$, Alice chooses at most $6M+1$ balls,
whose radii satisfy
\[
\sum_j\rho_{n,j}^s
\le C_\beta(\delta r_n/\beta)^s+\tfrac12(\beta r_n)^s
\le(\beta r_n)^s.
\]
Thus Alice's choices satisfy \eqref{eq:budget} and are legal.

Suppose $r_n\to0$ and the outcome $x_\infty$ avoids every ball
chosen by Alice; otherwise Alice wins. 
Then $x_\infty\ne z_J$ for all vertices. We show that
$x_\infty\in\widehat S$. For any vertex with $|J|<r_1$, let $n \geq 2$ be such that 
\[
r_n \leq |J| <  r_{n-1} \mbox{ and } r_n\ge\beta r_{n-1}>\beta|J| \geq \delta |J|.
\]
Consequently $r_n\le|J|<r_n/\beta$. If the marked ball $\overline{B}(z_J, \delta |J|) \cap B_n \neq \emptyset$, Alice has chosen it;  for otherwise $x_{\infty} \notin \overline{B}(z_J, \delta |J|)$ since $x_{\infty} \in B_n$. In either case,
\[
|x_\infty-z_J|\ge\delta|J|.
\]

Let $\mathcal F=\{J\in\T:|J|\ge r_1\}$. This family is finite:
a depth-$d$ vertex has length at most $q^d|I|$, so only finitely
many depths occur, and at each depth interior disjointness bounds
the number of such intervals by $|I|/r_1$. Assume that $\mathcal F\ne\varnothing$. Since  $x_\infty\ne z_J$ by Alice's choices and $|J| > 0$ for every $J\in\mathcal F$. Then we have that each ratio
$|x_\infty-z_J|/|J|$ is positive, and finiteness gives
\[
\eta:=\min_{J\in\mathcal F}\frac{|x_\infty-z_J|}{|J|}>0.
\]
If $\mathcal F=\varnothing$, set $\eta=\delta$.
Combining all these together, we obtain
\[
\inf_{J\in\T}\frac{|x_\infty-z_J|}{|J|}
\ge\min\{\delta,\eta\}>0.
\]
Hence $x_\infty\in\widehat S$.
This proves that $\widehat S$ is $0$-potential winning, and the absolute winning property follows from Proposition~\ref{prop:conversion}.

To see that $\widehat S\cap I^\circ\ne\varnothing$, play the
$\beta$-absolute game with $0<\beta<1/3$ and $B_1\subset I^\circ$.
Let Bob take $r_{n+1}=\beta r_n$ at every turn, which is possible
since $\beta<1/3$. Then $r_n\to0$, and Alice's winning strategy gives
\[
\bigcap_{n\ge1}B_n=\{x_\infty\}\subset\widehat S\cap I^\circ.
\]
\end{proof}

\begin{lemm}\label{lem:finite}
Let $\T$ be an interval tree with root interval $I$, and let
$L\ge1$ be an integer.
For each $J\in\T$, let $Z_J\subset J$ be a nonempty finite set
of points with $\card Z_J\le L$. Then the set
\[
\left\{x\in\R:
\inf_{J\in\T}\frac{\dist(x,Z_J)}{|J|}>0\right\}
\]
is absolute winning on $\R$, and its restriction to $I$ is absolute
winning.
\end{lemm}

\begin{proof}
For each $J\in\T$, choose $z_{J,1},\ldots,z_{J,L}\in J$ such that
\[
Z_J=\{z_{J,1},\ldots,z_{J,L}\},
\]
allowing repetitions. For each $\ell\in\{1,\ldots,L\}$, define
\[
S_\ell=
\left\{x\in\R:
\inf_{J\in\T}\frac{|x-z_{J,\ell}|}{|J|}>0\right\}.
\]
By Proposition~\ref{prop:tree}, $S_\ell$ is absolute
winning on $\R$. Hence $\bigcap_{\ell=1}^{L}S_\ell$ is absolute
winning by Proposition~\ref{prop:aw-properties}(2). Since
\[
\left\{x\in\R:
\inf_{J\in\T}\frac{\dist(x,Z_J)}{|J|}>0\right\}
=\bigcap_{\ell=1}^{L}S_\ell,
\]
this set is absolute winning on $\R$, and its intersection with $I$ is absolute winning on $I$.
\end{proof}

\section{The sequential return tree}\label{sec:lsv}

Fix $\aseq\in K^{\Nzero}$. We shall construct an interval tree
associated with $\aseq$, denoted by $\T_{\aseq}$.
We write $\T$ when the parameter sequence is fixed.
Define the countable set of preimages of the branch endpoints by
\begin{equation}\label{eq:boundary}
\mathcal B=
\bigcup_{n\ge0}(\Ta^n)^{-1}(\{0,1/2,1\}).
\end{equation}

Set $Y=[1/2,1]$. For $t\in\Nzero$ and integers $n \ge1$, define
\[
F_{t,n}=T_{a_{t+n-1}}\circ\cdots\circ T_{a_t}, \ 
F_{t,0}=\mathrm{id}.
\]
In particular, $F_{0,n}=\Ta^n$ for every $n\in\Nzero$.

\begin{lemm}\label{lem:entrance}
For every $t\in\Nzero$ and $u\in(0,1/2)$, there exists an integer $n\ge1$ such that
\[
F_{t,n}(u)\in Y.
\]
Consequently, the orbit of each $x\in[0,1]\setminus\mathcal B$
visits $Y^\circ$ infinitely often.
\end{lemm}
\begin{proof}
Fix $t\in\Nzero$ and $u\in(0,1/2)$, and set
\[
u_0=u, \  u_{k+1}=T_{a_{t+k}}(u_k) \ (k\ge0).
\]
Whenever $u_0,\ldots,u_k\in(0,1/2)$, we have
\[
u_{k+1}-u_k
=u_k(2u_k)^{a_{t+k}}
\ge u_k(2u_k)^{a_+}
\ge u(2u)^{a_+}=:d_u>0.
\]
Hence the orbit enters $Y$ in finitely many steps, with a bound
depending only on $u$ and $a_+$.

For $x\notin\mathcal B$, its forward orbit first visits $Y$
and then leaves it in finitely many steps, since the right branch
has slope 2. The preceding argument guarantees another
return, and repeating this argument gives infinitely many visits to $Y^\circ$,
since the orbit avoids $\{1/2,1\}$.
\end{proof}

We use the first return partition from \cite{BBD,BBR}.
The difference is that we take the closures of the partition intervals
when constructing the interval tree.

\subsection{Induced map}

For $t\in\Nzero$ and $u\in Y\setminus\{1/2\}$, define the {\it first return time} by
\[
r_t(u)=\min\{n\ge1:F_{t,n}(u)\in Y\}.
\]

For $t\in\Nzero$ and integers $\ell\ge1$, define
\begin{equation}\label{eq:return-itinerary}
V_{t,\ell}^{\circ}
:=\left(
T_{a_{t+\ell-1}}|_{(0,1/2)}\circ\cdots\circ
T_{a_{t+1}}|_{(0,1/2)}\circ T_{a_t}|_{Y^\circ}
\right)^{-1}(Y^\circ)
\end{equation}
and
\[
V_{t,\ell}:=\overline{V_{t,\ell}^{\circ}},
\]
where for $\ell=1$ the composition consists only of $T_{a_t}|_{Y^\circ}$.
For each $u\in V_{t,\ell}^{\circ}$,
$F_{t,k}(u)\in(0,1/2)$ for $1\le k<\ell$, and
$F_{t,\ell}(u)\in Y^\circ$. Hence $r_t(u)=\ell$.
For each fixed $t\in\Nzero$, the family
\[
\mathcal P_t=\{V_{t,\ell}^{\circ}:\ell\ge1\}
\]
is a partition of $Y$ modulo a countable set by Lemma~\ref{lem:entrance}.
We use the closed intervals
$V_{t,\ell}$ for the tree; these elements may have common endpoints.

Define the {\it induced map}
\[
R_t^{\circ}:\bigcup_{\ell\ge1}V_{t,\ell}^{\circ}\longrightarrow Y^\circ \ \mbox{ by } \  R_t^{\circ}|_{V_{t,\ell}^{\circ}}
=F_{t,\ell}|_{V_{t,\ell}^{\circ}}.
\]
Each branch is an increasing $C^1$ diffeomorphism onto $Y^\circ$.
Let $R_{t,\ell}:V_{t,\ell}\to Y$ be the unique $C^1$ extension of
$F_{t,\ell}|_{V_{t,\ell}^{\circ}}$. At an endpoint $b\in\partial V_{t,\ell}$, the extension value
$R_{t,\ell}(b)$ need not equal the actual iterate $F_{t,\ell}(b)$.

\subsection{The interval tree}

Take $I=[0,1]$ as the root. The vertices at depth one are the
closures of the first entry branches from $(0,1/2)$ into $Y$,
together with $Y$ itself. The entry time on $Y$ is 0.
For entry time $e\ge1$, the branch interior is
\[
\left(
T_{a_{e-1}}|_{(0,1/2)}\circ\cdots\circ
T_{a_0}|_{(0,1/2)}
\right)^{-1}(Y^\circ).
\]
For an interval $J$ at depth one with entry time $e$, set
$t(J):=e$ and let $H_J:J\to Y$ be the unique $C^1$ extension of
$\Ta^e|_{J^\circ}$. The map $H_J$ is an increasing diffeomorphism.
In particular,
\[
t(Y)=0, \  H_Y=\mathrm{id}_Y.
\]

Inductively, let $J$ be a nonroot vertex and put $t=t(J)$.
Pull the return partition $\mathcal P_t$
on $Y$ back to $J$ using $H_J$ and take closures.
For each $\ell\ge1$, define the child interval $J'$ of $J$,
its total iteration time $t(J')$, and its branch map $H_{J'}$ by
\begin{equation}\label{eq:children}
J'=H_J^{-1}(V_{t,\ell}), \ 
t(J')=t(J)+\ell, \ 
H_{J'}=R_{t,\ell}\circ H_J|_{J'}.
\end{equation}
Thus, for every nonroot vertex $J$, the map $H_J:J\to Y$ is the unique
$C^1$ extension of $\Ta^{t(J)}|_{J^\circ}$.

A vertex $J$ of depth $k\ge1$ specifies a first entry time $e=e(J)$ and
$k-1$ successive return times $\ell_1,\ldots,\ell_{k-1}$, with
\[
t(J)=e+\ell_1+\cdots+\ell_{k-1}.
\]
By construction, every child interval is contained in its parent,
and distinct children of the same vertex have disjoint interiors.
Thus the nesting and disjointness conditions in
Definition~\ref{def:tree} hold.

For $x\in[0,1]\setminus\mathcal B$, define
\[
s_0(x)=\min\{n\ge0:\Ta^n(x)\in Y\},
\]
and recursively,
\[
s_{i+1}(x)=\min\{n>s_i(x):\Ta^n(x)\in Y\}, \ i\ge0.
\]
These times are finite by Lemma~\ref{lem:entrance}.
Moreover, 
\[
s_{i+1}(x)
=s_i(x)+r_{s_i(x)}\bigl(\Ta^{s_i(x)}(x)\bigr), \ i\ge0.
\]
For each $k\ge1$, there is a unique $J_k(x)\in\T_k$
such that $x\in J_k(x)^\circ$, and
\[
t(J_k(x))=s_{k-1}(x), \ 
H_{J_k(x)}(x)=\Ta^{s_{k-1}(x)}(x).
\]

The tree $\T$ constructed as above is countable, and every vertex has finite depth.
By construction,
\[
\partial J\subset\mathcal B
\ \text{for every }J\in\T.
\]

\subsection{Distortion}

The following lemma applies the distortion estimate of
\cite{BQT} to the branch maps $H_J$.
For similar distortion estimates for small random perturbations
of intermittent circle maps, see \cite{SvS}.

\begin{lemm}[Uniform distortion]\label{lem:distortion}
There exists $D=D(K)\ge1$ such that for every $\aseq\in K^{\Nzero}$
and every nonroot vertex $J\in\T_{\aseq}$,
\begin{equation}\label{eq:distortion}
\frac{\sup_J H_J'}{\inf_J H_J'}\le D,
 \ 
\inf_J H_J'\ge\frac1{2D|J|},
 \  H_J(\inf J)=\frac12.
\end{equation}
Endpoint derivatives are understood as one-sided derivatives from interior.
\end{lemm}

\begin{proof}
The case $J=Y$ is trivial. Set
\[
\widetilde a_-=\min\{a_-/2,1/2\}, \ 
\widetilde K=[\widetilde a_-,a_+].
\]
Then $K\subset\widetilde K$, $0<\widetilde a_-<1$, and
$\widetilde a_-<a_+$, so \cite[Corollary~3.3]{BQT} applies to
$\aseq\in \widetilde K^{\Nzero}$.
The estimate is stated for branches with images contained in $[1/2,1)$.
By continuity, it also holds for their branch extensions with images
in $[1/2,1]$. Applying it to $H_J$ for $t(J)\ge1$ gives the first
inequality in \eqref{eq:distortion} with $D=D(K)\ge1$,
independent of $\aseq$ and $t(J)$.
Since $H_J:J\to Y$ is increasing and onto,
$H_J(\inf J)=1/2$ and
\[
\frac1{2|J|}
=\frac1{|J|}\int_J H_J'(x)\,dx
\le\sup_J H_J'\le D\inf_J H_J'.
\]
This completes the proof.
\end{proof}

\begin{lemm}\label{lem:child}
The interval tree $\T_{\aseq}$ satisfies \eqref{eq:contraction} with
$1/2\le q=1-(2D)^{-1}<1$.
\end{lemm}

\begin{proof}
Since $|I|=1$ and $q\ge1/2$, each interval $J$ at depth one
satisfies $|J|\le1/2\le q|I|$. Moreover, for each $t \in \Nzero$ and $\ell \geq 1$, $R_{t,\ell}\colon V_{t,\ell}\to Y$ satisfies $R_{t,\ell}^\prime\ge2$, we have $|V_{t,\ell}|\le|Y|/2$.

Now let $J'$ be a child of a nonroot vertex $J$, and put $t=t(J)$.
By \eqref{eq:children}, $H_J(J')=V_{t,\ell}$ for some $\ell\ge1$.
By Lemma~\ref{lem:distortion},
\[
\sup_Y(H_J^{-1})'\le D\inf_Y(H_J^{-1})'.
\]
Hence
\[
\frac{|J'|}{|J|}
=1-\frac{\int_{Y\setminus V_{t,\ell}}(H_J^{-1})'}
           {\int_Y(H_J^{-1})'}
\le1-\frac{|Y\setminus V_{t,\ell}|}{D|Y|}
\le1-\frac1{2D}=q.
\]
This finishes the proof.
\end{proof}

\section{Proof of Theorem A}\label{sec:avoidance}

Throughout this section, fix $\aseq\in K^{\Nzero}$ and write
$\T=\T_{\aseq}$.

\subsection{The neutral point}\label{sec:zero}
\begin{prop}\label{prop:zero}
The set $\Ea(0)$ is absolute winning on $I=[0,1]$.
\end{prop}

\begin{proof}
For each $J\in\T$, choose its left endpoint as the marked point
$z_J:=\inf J$. In particular, $z_I=0$.
By Proposition~\ref{prop:tree}, Lemma~\ref{lem:child}, and Proposition~\ref{prop:aw-properties}(3),
the set
\[
S_0=\left\{x\in[0,1]:\inf_{J\in\T}\frac{|x-z_J|}{|J|}>0\right\}
\setminus\mathcal B
\]
is absolute winning. Fix $x\in S_0$ and put
$\delta_x=\inf_{J\in\T}|x-z_J|/|J|>0$.

By Lemma~\ref{lem:entrance}, the orbit of $x$ visits $Y^\circ$ infinitely often.
For every $n\ge0$ with $\Ta^n(x)\in Y$, the construction in
Section~\ref{sec:lsv} gives a nonroot vertex $J$ such that
$x\in J^\circ$ and $t(J)=n$. Hence \eqref{eq:distortion} gives
\[
\Ta^n(x)-\frac12
=\int_{z_J}^{x}H_J'(u)\,du
\ge\frac{\delta_x}{2D}.
\]
After one time iterate, if $\Ta^{n+1}(x)\notin Y^\circ$, then
\[
\Ta^{n+1}(x)=2\Ta^n(x)-1\ge\delta_x/D.
\]
For otherwise, $\Ta^{n+1}(x) \geq 1/2$. On the other hand, whenever $\Ta^n(x)<1/2$, we have
$\Ta^{n+1}(x)\ge\Ta^n(x)$. Hence
\[
\inf_{n\ge0}\Ta^n(x)\ge
\min\{x,1/2,\delta_x/D\}>0.
\]
Thus $S_0\subset\Ea(0)$, so $\Ea(0)$ is absolute winning
by Proposition~\ref{prop:aw-properties}(1).
\end{proof}

\begin{lemm}\label{lem:gaps}
For each $x\in\Ea(0)\setminus\mathcal B$, there exists an integer
$N(x)\ge1$ such that the visiting times $s_k(x)$ defined in
Section~\ref{sec:lsv} satisfy
\[
s_{k+1}(x)-s_k(x)\le N(x)\ \text{for every }k\ge0.
\]
\end{lemm}
\begin{proof}
Let $\varepsilon=\min\{1/4,\inf_{n\ge0}\Ta^n(x)\}>0$.
Whenever $\Ta^n(x)<1/2$, we have
\[
\Ta^{n+1}(x)-\Ta^n(x)
\ge\varepsilon(2\varepsilon)^{a_+}.
\]
For each $k \geq 0$, if $\Ta^{s_k(x)+1}(x)\in Y$, then
$s_{k+1}(x)-s_k(x)=1$.
Otherwise, the preceding inequality gives
\[
s_{k+1}(x)-s_k(x)\le
1+\left\lceil\frac{1}{2\varepsilon(2\varepsilon)^{a_+}}\right\rceil
=:N(x).
\]
\end{proof}

\subsection{Estimates after return}\label{sec:offset}

Fix $j\in\Nzero$. For each nonroot vertex $J$,
put $t=t(J)$ and consider $F_{t,j}\circ H_J$.
Divide $J^\circ$ at the points $x$ satisfying
$F_{t,r}(H_J(x))=1/2$ for some $0\le r<j$.
This gives at most $2^j$ open branch intervals.
Let $V$ be the closure of any such interval, and let
$\phi_{J,j,V}:V\to[0,1]$ be the $C^1$ extension of
$(F_{t,j}\circ H_J)|_{V^\circ}$.

For $x\in V^\circ$,
\[
\phi_{J,j,V}(x)=F_{t,j}(H_J(x))=\Ta^{t+j}(x).
\]
The chain rule and Lemma~\ref{lem:distortion} give
\begin{equation}\label{eq:offsetderivative}
\begin{aligned}
\phi_{J,j,V}'(x)
&=H_J'(x)\prod_{r=0}^{j-1}
T_{a_{t+r}}'\bigl(\Ta^{t+r}(x)\bigr)\\
&\ge H_J'(x)\ge\frac1{2D|J|},
\end{aligned}
\end{equation}
since each factor in the product is at least $1$; for $j=0$ the
product equals $1$. By continuity of the branch derivative,
\eqref{eq:offsetderivative} also holds on $\partial V$.

For each $p\in[0,1]$ and each nonroot vertex $J$, define the finite set
\[
Z_J^{(j,p)}=
\bigcup_V\bigl(\partial V\cup\phi_{J,j,V}^{-1}(\{p\})\bigr)
\subset J,
\]
where the union is over these closed branch intervals. In particular, $\phi_{J,j,V}^{-1}(\{p\})$ could be empty set.
For the root interval put $Z_I^{(j,p)}=\{0,1\}$. 
Each $V$ contributes its two endpoints and at most one preimage
of $p$, since $\phi_{J,j,V}$ is strictly increasing. Hence for each $J\in\T$, 
\begin{equation}\label{eq:numbermarks}
1\le\card Z_J^{(j,p)}\le3\cdot2^j.
\end{equation}

\begin{prop}\label{prop:offset}
For every integer $j \in \Nzero$ and every $p \in [0, 1]$, the set
\[
S_{j,p}=
\left\{x\in[0,1]:
d_{j,p}(x):=
\inf_{J\in\T}
\frac{\dist(x,Z_J^{(j,p)})}{|J|}>0\right\}
\]
is absolute winning on $[0,1]$.
For every $x\in S_{j,p}\setminus\mathcal B$, the estimate
\begin{equation}\label{eq:offsetbound}
|\Ta^{t(J)+j}(x)-p|\ge\frac{d_{j,p}(x)}{2D}
\end{equation}
holds for every nonroot vertex $J$ with $x\in J^\circ$.
\end{prop}

\begin{proof}
The absolute winning result follows from
\eqref{eq:numbermarks} and Lemma~\ref{lem:finite}.
Fix $x\in S_{j,p}\setminus\mathcal B$ and a nonroot vertex $J$
with $x\in J^\circ$.
Since $x\notin\mathcal B$, the partition of $J^\circ$ constructed above
has a unique open branch interval containing $x$.
Let $V$ denote its closure and write $\phi=\phi_{J,j,V}$.

If $p\in\phi(V)$, set $z=\phi^{-1}(p)$.
If $p\notin\phi(V)$, choose $z\in\partial V$ such that
$|\phi(z)-p|=\dist(p,\phi(V))$.
In either case, $z\in Z_J^{(j,p)}$ and
$|\phi(x)-p|\ge|\phi(x)-\phi(z)|$.
Since $\phi(x)=\Ta^{t(J)+j}(x)$, the mean value theorem and
\eqref{eq:offsetderivative} give
\[
|\Ta^{t(J)+j}(x)-p|
\ge\frac{\dist(x,Z_J^{(j,p)})}{2D|J|}
\ge\frac{d_{j,p}(x)}{2D}.
\]
\end{proof}

\subsection{Proof of Theorem~\ref{thm:main}}\label{sec:proof}
\begin{proof}[Proof of Theorem~\ref{thm:main}]
Fix $p\in I$. Define
\begin{equation}\label{eq:exception}
\mathcal C_p=
\bigcup_{n\ge0}(\Ta^n)^{-1}(\{p,0,1/2,1\}).
\end{equation}
Since each $T_a$ has two injective branches, every point has
at most $2^n$ preimages under $\Ta^n$.
Hence $\mathcal C_p$ is countable and contains $\mathcal B$.
Consider
\begin{equation}\label{eq:winningcore}
S=
\left(\Ea(0)\cap
\bigcap_{j=0}^{\infty}S_{j,p}\right)
\setminus\mathcal C_p.
\end{equation}
By Proposition~\ref{prop:zero}, Proposition~\ref{prop:offset},
and Proposition~\ref{prop:aw-properties}(2)--(3),
$S$ is absolute winning on $I$.

Choose any $x\in S$, let $N(x)<\infty$ be the integer given by Lemma~\ref{lem:gaps}. From the definition of $S$, we have $d_{j,p}(x)>0$ for each $j\ge0$. Then let
\begin{equation}\label{eq:finitemin}
\eta_p(x)=\frac1{2D}
\min_{0\le j\le N(x)}d_{j,p}(x)>0.
\end{equation}
For $n\ge s_0(x)$, choose $k\ge0$ such that
$s_k(x)\le n<s_{k+1}(x)$, and put $j=n-s_k(x)$.
Then $0\le j<N(x)$.
Let $J=J_{k+1}(x)$ be the unique vertex at depth $k+1$
whose interior contains $x$. Then $t(J)=s_k(x)$.
Applying Proposition~\ref{prop:offset} with this vertex $J$ and
$j=n-s_k(x)$, so that $t(J)+j=n$, gives
\[
|\Ta^n(x)-p|
\ge\frac{d_{j,p}(x)}{2D}
\ge\eta_p(x).
\]

It remains to consider $0\le n<s_0(x)$. If $s_0(x)>0$, none of the finitely many values
$\Ta^n(x)$, $0\le n<s_0(x)$, equals $p$, since
$x\notin\mathcal C_p$. Hence
\[
\tilde\eta_p(x)
=\min_{0\le n<s_0(x)}|\Ta^n(x)-p|>0.
\]
For $s_0(x)=0$, set $\tilde\eta_p(x)=1$.
Combining these estimates together give
\[
\inf_{n\ge0}|\Ta^n(x)-p|
\ge\min\{\eta_p(x),\tilde\eta_p(x)\}>0.
\]
Thus $S\subset\Ea(p)$, so $\Ea(p)$ is absolute winning
by Proposition~\ref{prop:aw-properties}(1).

For every nonempty countable set $P\subset I$, the identity
\[
\Ea(P)=\bigcap_{p\in P}\Ea(p)
\]
and Proposition~\ref{prop:aw-properties}(2) show that $\Ea(P)$
is absolute winning. If $P=\varnothing$, then $\Ea(P)=I$,
and the conclusion is trivial.
\end{proof}

\subsection*{AI Declaration}
The authors used ChatGPT and Codex (OpenAI; models GPT-6 Astra and
GPT-6.1 Sol) to assist with language editing, notation checks and reference searches.
ChatGPT also suggested organizing the inverse branches into an interval tree, which simplified the presentation of the proof.
All other statements and proofs in the paper were developed
by the authors, who take full responsibility for the contents of the paper.

\begingroup
\raggedright

\endgroup

\bigskip
\noindent\begin{minipage}{\textwidth}
\small
Haoyang Ji\\
School of Mathematics and Statistics, Zhengzhou University\\
Zhengzhou 450001, China\\
E-mail: \href{mailto:jihymath@zzu.edu.cn}{jihymath@zzu.edu.cn}
\end{minipage}

\end{document}